\documentclass[authoryear,11pt]{article}
\usepackage{ifpdf}
\usepackage[latin1]{inputenc}
\usepackage{tikz}
\usepackage{verbatim}
\usetikzlibrary{arrows,shapes}

\usepackage{graphicx,natbib,amssymb,lineno}
\usepackage{times}
\usepackage{color}

\def \bX {\textbf{X}}
\def \cind {\,\bot\hspace{-.50em}\bot\,}
\def \ncind {\,\not\!\!\!\bot\hspace{-.50em}\bot\,}

\def \r {{\rm I\hspace{-.10em}R}}

\newenvironment{pf}[1][Proof]{\noindent\textbf{#1.} }{\ \rule{0.5em}{0.5em}\\}

\newtheorem{theorem}{Theorem}
\newtheorem{definition}{Definition}
\newtheorem{lemma}[theorem]{Lemma}
\newtheorem{example}{Example}

\newtheorem{proposition}[theorem]{Proposition}
\newenvironment{rk}[1][Remark]{\noindent\textbf{#1.} }{\ \\}

\title{Gaussian Covariance faithful Markov Trees }

\author{Dhafer Malouche \\
{\small Ecole Supérieure de la Statistique}\\
{\small et de l'Analyse de l'Information, Tunisia. }\\
{\small \texttt{dhafer.malouche@essai.rnu.tn}}\\
\\
and \\
\\
Bala Rajaratnam \\
{\small Standford University, USA. }\\
{\small\texttt{brajarat@stanford.edu}}}

\date{}

\begin{document}

\maketitle
\begin{abstract}

 A covariance graph is an undirected graph associated with a multivariate probability distribution of a given random vector  where each vertex represents each of the different components of the random vector and  where the absence of an edge between any pair of variables implies marginal independence between these two variables. Covariance graph models have recently received much attention in the literature and constitute a sub-family of graphical models. Though they are conceptually simple to understand, they are considerably more difficult to analyze. Under some suitable assumption on the probability distribution, covariance graph models can also be used to represent more complex conditional independence relationships between subsets of variables. When the covariance graph captures or reflects  all the conditional independence statements present in the probability distribution the latter is said to be faithful to its covariance graph - though no such prior guarantee exists. Despite the increasingly widespread use of these two types of graphical models, to date no deep probabilistic analysis of this class of models, in terms of the faithfulness assumption, is available. Such an analysis is crucial in understanding the ability of the graph, a discrete object, to fully capture the salient features of the probability distribution it aims to describe. In this paper we demonstrate that multivariate Gaussian distributions that have trees as covariance graphs are necessarily faithful. The method of proof is original as it uses an entirely new approach and in the process yields a technique that is novel to the field of graphical models.
\end{abstract}

\section{Introduction}
Markov random fields or graphical models are widely used to represent conditional independences in a given multivariate probability distribution (see \cite{KunschGemanKehagias}, \cite{JiSeymour}, \cite{Spitzer}, \cite{KindermannSnell}, \cite{Lauritzen1996} to name just a few). Many different types of Markov Random fields or graphical models have been studied in the literature. For example,  directed acyclic  graphs or DAGs are  commonly referred to as ``Bayesian networks" (see \cite{Pearl1988}). When the graph is undirected and when such graphs are constructed using marginal independence relationships between pairs of random variables in a given random vector these graphical models are called ``covariance graph" models (see \cite{Cox1993}, \cite{Cox1996}, \cite{Kauermann1996}, \cite{Malouche2009} and \cite{Khare2009}). Covariance graph models are commonly represented by graphs with exclusively bi-directed or dashed edges (see \cite{Kauermann1996}). This representation is used in order to distinguish them from the traditional and widely used concentration graph models. Concentration graphs encode conditional independence between pairs of variables given the remaining ones. Formally, if we  consider a random vector  $\mathbf{X}=(X_v, v\in V)'$ with a probability distribution $P$  where $V$ is a finite set representing the random variables in $\mathbf{X}$. The concentration graph associated with $P$ is an undirected graph $G=(V,E)$ where
\begin{itemize}
\item $V$ is the set of vertices.
\item Each vertex represents one variable in $\mathbf{X}$.
\item $E$ is the set of edges (between the verices in $V$) constructed using the pairwise rule : for pair $(u,v)\in V\times V$, $u\not=v$
\begin{equation}
\label{defCon}
(u,v)\not\in E\;\iff\; X_u\cind X_v\mid \mathbf{X}_{V\setminus\{u,v\}}
\end{equation}
where $\mathbf{X}_{V\setminus\{u,v\}}:=(X_w,\, w\not=u\mbox{ and }w\not= v)'$.
\end{itemize}
Note that $(u,v)\not\in E$ means that the vertices $u$ and $v$ are not adjacent in $G$.

An undirected graph $G_0=(V,E_0)$ is called the covariance graph associated with the probability distribution $P$ if the set of edges $E_0$ is constructed as follows
\begin{equation}
\label{defCov}
(u,v)\not\in E\;\iff\; X_u\cind X_v
\end{equation}
The subscript zero is invoked for covariance graphs  (i.e., $G_0$ vs $G$) as  the definition of covariance graphs does not involve  conditional independences.

Both concentration and covariance graphs are not only used to encode pairwise relationships between pairs of variables in the random vector $\mathbf{X}$, but as we will see below, these graphs can be used to  encode conditional independences  that exist between subsets of variables of $\mathbf{X}$. First we introduce some definitions:

The multivariate distribution $P$ is said to satisfy the ``intersection property" if for any subsets $A$, $B$ $C$ and $D$ of $V$ which are pairwise disjoint,
\begin{equation}\label{inter}
\left\{\begin{array}{lcl}
     \mathbf{X}_A\cind \mathbf{X}_B\mid \mathbf{X}_{C\cup D} & & \\
\mbox{and }  &\mbox{ then } &   \mathbf{X}_A\cind \mathbf{X}_{B\cup C}\mid \mathbf{X}_{ D}\\
\mathbf{X}_A\cind \mathbf{X}_C\mid \mathbf{X}_{B\cup D} & & \\
\end{array}\right.\end{equation}

We will call the \textit{intersection} property (see \cite{Lauritzen1996}) in (\ref{inter}) above the \textit{concentration intersection} property in this paper in order to differentiate it from another  property that is satisfied by $P$ when studying  covariance graph models.

Let $P$ satisfy the concentration intersection property. Then for any triplet $(A,B,S)$ of subsets of $V$ pairwise disjoint, if $S$ separates\footnote{We say that $S$ separates $A$ and $B$ if any path connecting $A$ and $B$ in $G$ intersects $S$, i.e., $A\bot_G B\mid S$, and is not to be confused with stochastic independence which is denoted by $\cind$  as compared to $\bot_G$.} $A$ and $B$ in the concentration graph $G$ associated with $P$ then the random vector $\mathbf{X}_A=(X_v,\, v\in A)'$ is independent of $\mathbf{X}_B=(X_v,\, v\in B)'$ given $\mathbf{X}_S=(X_v,\, v\in S)'$. This latter property is called \textit{concentration global  Markov} property and is formally defined as,
\begin{equation}\label{conGmarkov}
A\bot_G B\mid S\; \Rightarrow\;\mathbf{X}_A\cind \mathbf{X}_B\mid \mathbf{X}_S.
\end{equation}
\cite{Kauermann1996} and  \cite{Banerjee2003} show that if $P$ satisfies the following property :
for any triplet $(A,B,S)$ of subsets of $V$ pairwise disjoint,
 \begin{equation}\label{interCov}
\mbox{ if } \mathbf{X}_A\cind \mathbf{X}_B\mbox{ and }\mathbf{X}_A\cind \mathbf{X}_C\;\mbox{ then }\mathbf{X}_A\cind \mathbf{X}_{B\cup C},
\end{equation}
then for any triplet $(A,B,S)$ of subsets of $V$ pairwise disjoint, if $V\setminus(A\cup B\cup S)$ separates $A$ and $B$ in the covariance graph $G_0$ associated with $P$ then $\mathbf{X}_A\cind \mathbf{X}_B\mid \mathbf{X}_S$. This latter property is called the \textit{covariance global Markov  property} and can be written formally as follows
 \begin{equation}A\bot_{G_0} B\mid V\setminus(A\cup B\cup S)\;\Rightarrow\;\mathbf{X}_A\cind \mathbf{X}_B\mid \mathbf{X}_S.\label{covGmarkov}
 \end{equation}
In parallel to  the concentration graph case,  property (\ref{interCov}) will be called the \textit{covariance intersection} property.

Even if $P$ satisfies both intersection properties, the covariance and concentration graphs may not be able to capture or reflect all the conditional independences present in the distribution, i.e., there may exist one or more conditional independences present in the probability distribution that does not correspond to any separation statement in either $G$ or $G_0$. Equivalently, a lack of a separation statement in the graph does not necessarily imply conditional independences. On the contrary case when no other conditional independence exist in $P$ except the ones encoded by the graph, we classify $P$ as a \textit{faithful} probability distribution to its graphical model. More precisely we say that $P$ is  \textit{concentration faithful} to its concentration graph if for any   triplet $(A,B,S)$ of subsets of $V$ pairwise disjoint, the following statement holds :
\begin{equation}\label{conFaith}
S\mbox{ separates }A\mbox{ and }B\;\iff\; \mathbf{X}_A\cind \mathbf{X}_B\mid \mathbf{X}_S.
\end{equation}
Similarly, $P$ is said to be \textit{covariance faithful} to its covariance graph $G_0$ if for any   triplet $(A,B,S)$ of subsets of $V$ pairwise disjoint, the following statement holds :
\begin{equation}\label{covFaith}
V\setminus(A\cup B\cup S)\mbox{ separates }A\mbox{ and }B\;\iff\; \mathbf{X}_A\cind \mathbf{X}_B\mid \mathbf{X}_S.
\end{equation}
A natural question of both theoretical and  applied interest in probability theory is to understand the implications of the \textit{faithfulness} assumption. This assumption is fundamental since it yields a bijection  between the probability distribution  $P$ and the graph $G$ in terms of  the independences that are present in the distribution. In this paper we show that when $P$ is a multivariate Gaussian distribution whose covariance graph are trees are necessarily covariance faithful, i.e., these probability distributions satisfy property (\ref{covFaith}), i.e., the associated covariance graph $G$ is fully able to capture all the conditional independences present in the multivariate distribution $P$. This result can be considered as a dual of a previous probabilistic  result proved by \cite{Becker2005} for concentration graphs that demonstrates  that Gaussian distributions having concentration trees, i.e., the concentration graph is a tree are necessarily \textit{concentration faithful} to its concentration graph (implying property (\ref{conFaith}) is satisfied). This  result was proved by showing that  Gaussian distributions satisfy an additional intersection property. The approach in the proof of the main  result of this paper is vastly different from the one used for concentration graphs  by \cite{Becker2005}.

The outline of this paper is follows. Section 2 presents graph theory preliminaries. Section 3 gives a brief overview of covariance and concentration graphs associated with multivariate Gaussian distributions. Furthermore, an easier way to encode conditional independence using covariance graphs is given in Section 3. The prove of the main result of this paper is given in Section 4. Section 5 concludes by summarizing the results in the paper and the implications thereof.

\section{Graph theory preliminaries}
This section introduces notation and terminology that is required in subsequent sections. An undirected graph $G=(V,E)$ consists of two sets $V$ and $E$, with $V$ representing the set of vertices, and $E\subseteq (V\times V)\setminus\{(u,u),\, u\in V\}$ the set of edges satisfying :
$$ \forall \;(u,v)\in E\,\iff \,(v,u)\in E
$$
For $u,\,v\in V$, we write $u\sim _{G}v$ when $(u,v)\in E$ and we
say that $u$ and $v$ are \textit{adjacent} in $G$.

\begin{definition}
A \textit{path} connecting  two  distinct vertices $u$ and $v$ in $G$ is a sequence of distinct vertices
    $\left(u_0,u_1,\ldots,u_n)\right)$ where $u_0=u$ and $u_n=v$ where for every  $i=0,\ldots,n-1$, $u_i\sim_G u_{i+1}.$
\end{definition}

Such a path will be denoted $p=p(u,v,G)$ and we  say that $p(u,v,G)$ connects $u$ and $v$ or alternatively $u$ and $v$ are \textit{connected} by $p(u,v,G)$. Its \textit{length}, denoted by $|p(u,v,G)|$,  is defined as the number of edges connecting  the vertices of $p$. So, in this case $|p(u,v,G)|=n$. We also denote by $\mathcal{P}(u,v,G)$ the set of paths between $u$ and $v$.

Trees are a particular class of graphs that are studied in this paper. This class of graphs are formally defined below.

\begin{definition}
Let $G=(V,E)$ be an undirected graph. The graph $G$ is called a tree if any pair of vertices $(u,v)$ in $G$ are connected  by exactly one path, i.e., $|\mathcal{P}(u,v,G)|=1 \; \; \forall \; u,v \in V$.
\end{definition}

  A   subgraph  of $G$  \textit{induced} by a subset $U\subseteq
V$ is   denoted by $G_U=(U,E_U)$, $U\subseteq V$ and $E_U=E\cap
(U\times U) $.

\begin{definition}
 A \textit{connected component} of a graph $G$ is the largest subgraph $G_U=(U,E_U)$ of $G$ such that each pair of vertices can be connected by at least one path in $G_U$.
\end{definition}

We now state a Lemma needed in the proof of the main result of this paper.

\begin{lemma}\label{TreeCon}Let $G=(V,E)$ be an undirected graph. If $G$ is a tree, any subgraph of $G$ induced by a subset of $V$ is a union of connected components, each of which are trees (or what we shall refer to as a ``union of tree connected components").
\end{lemma}
\begin{pf}
Consider $U\subset V$, the induced graph $G_U$ and  a pair of vertices $(u,v)\in U\times U$. Let us assume to the contrary that $u$ and $v$ are connected by two distinct paths $p_1$ and $p_2$ in $G_U$ (i.e., $G_U$ is not a tree). As the set of edges $E_U$ of the graph $G_U$ is included in the set of edges $E$ of $G$, i.e., $E_U=E\cap (U\times U)\subseteq E$, then $p_1$ and $p_2$ are also paths in $G$. Hence $u$ and $v$ are vertices in $G$ which are connected by two distinct  paths, i.e., $p_1$ and $p_2$. This of course yields a contradiction with the fact that $G$ is a tree. Thus any pair of vertices in $G_U$ are connected by at most one path and, hence $G_U$ is a  union of connected components, each of which are trees  (or a ``union of tree connected components").
\end{pf}

\begin{definition}
 For a connected graph, a \textit{separator} is a subset $S$ of $V$
such that there exists a pair of non-adjacent vertices $u$ and $v$
such that $u,$ $v\not\in S$ and
  \begin{equation}\label{ugraphs:eq2}
  \forall p\in \mathcal{P}(u,v,G),\;\; p\cap S\not=\emptyset
  \end{equation}
\end{definition}

If $S$ is a separator then it is easily verified  that every $S'\supseteq S$ such that $S'\subseteq V\setminus\{u,v\}$ is also a separator. We are thus lead to the notion of a minimal separator.

\begin{definition}
  The separator $S$ is defined to be  a \textit{minimal separator} between two non-adjacent vertices $u$ and $v$ if for any $w\in S$, the subsets $S\setminus\{w\}$ is not a separator of $u$ and $v$.
\end{definition}

Note that in the case where $G$ contains more than two connected components and if $u$ and $v$ belong to different connected components the empty set is the only possible separator of $u$ and $v$. Finally, let $A$, $B$ and $S$ be pairwise disjoint subsets of $V$. We say that $S$ \textit{separates} $A$ and $B$ if for any pair of vertices
$(u,v)\in A\times B$, any path connecting $u$ and $v$ intersects $S$. In the case where $A$ and $B$ belong to different connected components of $G$ the subset $S$ can be empty because the set of paths between any pair of vertices $(u,v)\in A\times B$ is empty.

\section{Gaussian Concentration and Covariance Graphs }
In this section we present a brief overview of concentration and covariance graphs in the case when the probability distribution $P$ is multivariate Gaussian. Such graphical models are commonly referred to as Gaussian covariance or Gaussian concentration graph models.

\subsection{Gaussian concentration graph models}
Consider a probability space with triplet $(\Omega,{\cal F},\mathbb{P})$ and let $\mathbf{X}\,:\, \Omega\rightarrow \mathbb{R}^{|V|}$ be a random vector where $\mathbf{X}=(X_v,\,v\in V)'$ and $P$ represents the induced measure of $\mathbb{P}$ by $\mathbf{X}$. If $\mathbf{X}$ follows a Gaussian distribution then it has  the following density function with respect to Lebesgue measure :
\begin{equation}\label{GaussDens}
 f(\mathbf{x})=\frac{1}{(2\pi)^{|V|/2}|\Sigma|^{1/2}} \,
  \exp\left(-\frac{1}{2}(\mathbf{x}-\mathbf{\mu})'\Sigma^{-1}(\mathbf{x}-\mathbf{\mu})'\right),
  \end{equation}
 where $\mathbf{x}=(x_u,\, u\in V)'\in \r^{|V|}$, $\mathbf{\mu}\in \r^{|V|}$ is the mean vector  and $\Sigma=(\sigma_{uv})\in \mathcal{P}^+$ is the  covariance matrix with $\mathcal{P}^+$ denoting the cone of symmetric positive definite matrices. Without loss of generality we will assume that $\mathbf{\mu}=\mathbf{0}$. As any Gaussian distribution with $\mathbf{\mu}=\mathbf{0}$ is completely determined by  its covariance matrix $\Sigma$, this set of multivariate Gaussian distributions can therefore be identified by the set of symmetric positive definite matrices.   Gaussian distributions can also be parameterized by the inverse of the covariance matrix $\Sigma$ denoted by $K=\Sigma^{-1}=(k_{uv})$. The matrix $K$ is called the \textit{precision} or \textit{concentration} matrix. It is well known (see \cite{Lauritzen1996}) that for any pair of variables $(X_u,X_v)$, where $u\not=v$
$$X_u\cind X_v\mid \mathbf{X}_{V\setminus\{u,v\}}\;\iff\; k_{uv}=0.$$
Hence  the concentration graph $G=(V,E)$ can be constructed simply using the precision matrix $K$ and the following  rule
$$(u,v)\not\in E\;\iff\; k_{uv}=0.$$
Furthermore it can be easily deduced from a classical result in \cite{Hammersly1971}, that is
reproved in \cite{Lauritzen1996}, that any multivariate random vector with a positive density necessarily satisfies
  the concentration intersection property (\ref{inter}). Hence for Gaussian concentration graph models the pairwise Markov property in (\ref{defCon}) is equivalent to the concentration global Markov property in (\ref{conGmarkov}).

\subsection{Gaussian covariance graph models}

As  seen earlier in (\ref{defCov}) covariance graphs are constructed using pairwise marginal independence  relationships. It is also well known that for  multivariate Gaussian distributions :
$$X_u\cind X_v\;\iff\, \sigma_{uv}=0.$$

Hence in the Gaussian case the covariance graph $G_0=(V,E_0)$ can be constructed using the following  rule :
$$(u,v)\not\in E_0\;\iff\; \sigma_{uv}=0.$$
It is also easily seen that Gaussian distributions satisfy the covariance intersection property defined in (\ref{interCov}).
Hence Gaussian covariance graphs can also encode conditional independences according to the following rule :
for any triplet $(A,B,S)$ of subsets of $V$ pairwise disjoint, if $V\setminus(A\cup B\cup S)$ separates $A$ and $B$ in the covariance
graph $G_0$  then $\mathbf{X}_A\cind \mathbf{X}_B\mid \mathbf{X}_S$. We now show (see proposition \ref{EquivGlob} below) that
 there is a simple way to read conditional independence statements from the covariance graph. This result holds
true for any probability distribution that satisfy the covariance intersection property given in (\ref{interCov}).

\begin{proposition}\label{EquivGlob}
Let $\bX_V=(X_v,\, v \in V)'$ be a random vector with probability
distribution $P$ satisfying the covariance intersection property in (\ref{interCov}) and let $G_0=(V,E_0)$ be the covariance  graph associated with $P$. Then
the following statements are equivalent,
\begin{itemize}
   \item[i.]  for any pairwise disjoint subsets $A$, $B$ and $S$ of $V$ : if $V\setminus(A\cup B\cup
S)$ separates $A$ and  $B$ in $G_0$ then  $\bX_A\cind \bX_B\mid
\bX_{S}$

 \item[ii.]  for any pairwise disjoint subsets $A$, $B$ and $S$ of $V$ : if $S$ separates $A$ and  $B$ in $G_0$ then  $\bX_A\cind \bX_B\mid
\bX_{V\setminus(A\cup B\cup S)}$

\end{itemize}
\end{proposition}

\begin{pf}
Let us first assume that (i) is satisfied and let us prove (ii).

Let $A$, $B$ and $S$ be three pairwise disjoint
 subsets of $V$ such that  $S$ separates $A$ and $B$ in $G_0$. Note that
 we can write $S$ as follows:
 $$S=V\setminus(V\setminus(A\cup B\cup S))\cup A \cup B)$$
 Since $(V\setminus(A\cup B\cup S)\cup A \cup B= V\setminus S$ and
 $V\setminus(V\setminus S)=S$.

 By hypothesis $S$ separates $A$ and $B$ in $G_0$. Let $S'=V\setminus(A\cup B\cup S)$ and since $S=V\setminus(S'\cup A\cup B)$ we can apply property (i) to the triplet $(A,B,S')$. Hence  $\bX_A\cind \bX_B\mid
 \bX_{S'}$. Hence $\bX_A\cind \bX_B\mid \bX_{V\setminus(S\cup A\cup B)}$ since $S':=V\setminus(S\cup A\cup
 B)$. We have therefore proved that  if $S$ separates $A$ and $B$ in $G_0$, then $\bX_A\cind \bX_B\mid
 \bX_{V\setminus(S\cup A\cup B)}$.

 Assume now that property (ii) is satisfied and let $A$, $B$ and $S$ be three pairwise disjoint
 subsets of $V$ such that   $V\setminus(S\cup A\cup B)$ separates
 $A$ and $B$ in $G_0$. Let us denote by $S'=V\setminus(S\cup A\cup B)$ which is a subset separating $A$ and $B$ in $G_0$.
 Since (ii) is satisfied, we deduce that $\bX_A\cind \bX_B\mid
 \bX_{V\setminus(A\cup B\cup S')}$. However
 $$V\setminus(A\cup B\cup S')=V\setminus((V\setminus(A\cup B\cup S))\cup A\cup
 B)=S$$
 Hence we conclude that $V\setminus(A\cup B\cup S)$ separates
 $A$ and $B$ in $G_0$ implies that  $\bX_A\cind \bX_B\mid
 \bX_{S}$. Thus  property (i) is satisfied.
\end{pf}

 Proposition \ref{EquivGlob} can be used to formulate  an equivalent definition of the covariance faithfulness property.

\begin{definition}
Let $\bX_V=(X_v,\, v \in V)'$ be a random vector with probability distribution $P$ satisfying the covariance intersection property in (\ref{interCov}) and let $G_0=(V,E_0)$ be the covariance  graph associated with $P$. We say that $P$ is covariance  faithful to $G_0$ if  for any pairwise disjoint subsets $A$, $B$ and $S$ of $V$ the following condition is satisfied $$S\mbox{ separates }A\mbox{ and }B\;\iff\; \bX_A\cind \bX_B\mid \bX_{V\setminus(A\cup B\cup S)}$$
\end{definition}
The above reformulation of the covariance faithfulness property is an important ingredient in the proofs in the next section.

\section{Gaussian Covariance faithful trees}

We now proceed to study  the faithfulness assumption in the context of multivariate Gaussian distributions and when the associated covariance graphs are trees.

The main result of this paper, presented in Theorem \ref{faithG0}, proves that multivariate Gaussian probability
distributions having tree covariance graphs are necessarily faithful  to their covariance graphs. The analogous  result for concentration graphs was demonstrated by \cite{Becker2005} where the authors proved  that Gaussian distributions having tree concentration graphs are necessarily faithful  to these graphs. We now formally state  Theorem \ref{faithG0}. The proof follows shortly after a series of lemmas/theorem(s) and an illustrative example.

\begin{theorem}\label{faithG0}
Let $\bX_V=(X_v,\, v\in V)'$ be a  random vector with Gaussian
distribution
 $P=\mathcal{N}_{|V|}(\mu,\Sigma^{-1})$. Let $G_0=(V,E_0)$ be  the covariance graph associated with $P$.
 If $G_0$ is a tree or more generally a union of connected components each of which are trees (or a union of ``tree connected components"),  then $P$ is $g_0-$faithful to $G_0$.
\end{theorem}

The proof of Theorem \ref{faithG0} requires among others a result proved by \cite{Jones2005}. This result gives a method
that can be used  to compute the covariance matrix $\Sigma$ from the precision matrix $K$ using the paths in the concentration graph $G$. The result can also be easily extended  to show that the precision matrix $K$ can be computed from the covariance matrix $\Sigma$ using the paths in the covariance graph $G_0$. We now state the result by \cite{Jones2005}.

\begin{theorem}\cite{Jones2005}.
\label{kgraphs:theo1}

 Let $\bX_V=(X_v,\, v\in V)'$ be a random vector with Gaussian distribution $P=\mathcal{N}_{|V|}(\mu,\Sigma)$ where $\Sigma$ and $K=\Sigma^{-1}$ are positive definite matrices. Let $G=(V,E)$ and $G_0=(V,E_0)$ denote respectively the concentration and covariance graph associated with the probability distribution of $\bX_V$.

For all $(u,v)$ in $V\times V$
$$k_{uv}=\displaystyle\sum_{p\in \mathcal{P}(u,v,G_0)}(-1)^{|p|+1}
|\sigma|_p\,\frac{|\Sigma\setminus p|}{|\Sigma|}$$ and
$$\sigma_{uv}=\displaystyle\sum_{p\in \mathcal{P}(u,v,G)}(-1)^{|p|+1} |k|_p\frac{|K\setminus
p|}{|K|}$$

 where, if $p=(u_0,\ldots,u_n),$
  $$|\sigma|_p=\sigma_{u_0u_1}\sigma_{u_1u_2}\ldots\sigma_{u_{n-1}u_n},\;\;
  |k|_p=k_{u_0u_1}k_{u_1u_2}\ldots k_{u_{n-1}u_n},$$
 $K\setminus
p=\left(k_{uv},\, (u,v)\in (V\setminus p)\times (V\setminus
p)\right)$ and $\Sigma\setminus p=\left(\sigma_{uv},\, (u,v)\in
(V\setminus p) \times (V\setminus p)\right)$ denote respectively $K$ and $\Sigma$
with rows and columns corresponding to variables in path $p$
omitted. The determinant of a zero-dimensional matrix is defined to
be $1$.

\end{theorem}

The proof of our main theorem (Theorem \ref{faithG0}) also requires the results proved in  the lemma  below.

\begin{lemma} \label{lemComplete}
 Let $\bX_V=(X_v,\, v\in V)'$ be a random vector with Gaussian
distribution
 $P=\mathcal{N}_{|V|}(\mu,K=\Sigma^{-1})$. Let $G_0=(V,E_0)$ and
 $G=(V,E)$ denote respectively  the covariance and  concentration
 graphs
 associated with $P$, then
\begin{itemize}
  \item[i.] $G$ and $G_0$ have the same connected components

  \item[ii.] If a given connected component in $G_0$ is a tree then the
  corresponding connected component in $G$ is complete and
  vice-versa.

\end{itemize}
\end{lemma}

\begin{pf}

\begin{itemize}
\item[] Proof of (i).

The fact that $G_0$ and $G$ have the same connected components can be  deduced from the matrix structure of the covariance and the precision matrix. The connected components of $G_0$ correspond to block diagonal matrices in $\Sigma$. Since $K=\Sigma^{-1}$, then by properties of inverting partitioned matrices, $K$  also has  the same  block diagonal matrices as $\Sigma$ in terms of the variables that constitute these  matrices. These blocks corresponds to distinct components in $G$ and $G_0$. Hence  both matrices have the same connected components.

\medskip

\item[] Proof of (ii).

 Let us assume now that the covariance graph $G_0$ is a tree, hence it is a connected graph with only one connected component. We shall prove that the concentration graph $G$ is complete by using Theorem \ref{kgraphs:theo1} by \cite{Jones2005} and  computing any coefficient $k_{uv}$ ($u\not =v$). Since $G_0$ is a tree,  there exists exactly one path between between any two vertices $u$ and $v$. We shall denote this path as $p=(u_0=u,\ldots,u_n=v)$. Then by Theorem \ref{kgraphs:theo1}
\begin{equation}
\label{faithGauss1} k_{uv} =(-1)^{n+1}\sigma_{u_0u_1}\ldots
\sigma_{u_{n-1}u_n}\displaystyle\frac{\left|\Sigma\setminus
p\right|} {\left|\Sigma\right|}\end{equation} First note that the
determinant of the matrices in (\ref{faithGauss1}) are all positive
since principal minors of positive definite matrices are positive.
Second since we are considering a path in $G_0$,  $\sigma_{u_{i-1}u_i}\not=0$, $\forall \; i=1,\ldots,n$. Using these two facts we deduce from (\ref{faithGauss1}) that  $k_{uv}\not=0$ for all $(u,v)\in E$. Hence $u$ and
$v$ are adjacent in $G$ for all $(u,v)\in E$. The concentration graph $G$ is  therefore complete. The proof that when $G$ is assumed to be a tree implying that $G_0$ is complete follows similarly.

\end{itemize}
\end{pf}

\begin{rk} We further note that Theorem \ref{kgraphs:theo1} is also directly useful in deducing the completeness of the concentration graph by using the covariance graph in other settings. As a concrete example consider the case when $G_0$ is a cycle with an even number of edges s.t. $|V|=2k$ for some  odd integer $k$, and assume that all the coefficients in the covariance matrix $\Sigma$ of $\bX_V$ are positive. Hence a given pair of vertices $(u,v)$ in $G_0$ are connected by two paths which are both of odd length. Let us  denote these paths as $p_1$ and $p_2$. Using Theorem \ref{kgraphs:theo1},
it is easily deduced that
$$k_{uv}=\sigma_{|p_1|}\frac{|\Sigma\setminus p_1|}{|\Sigma|}+\sigma_{|p_2|}\frac{|\Sigma\setminus p_2|}{|\Sigma|}$$
Here $|\sigma_{p_1}|$ and $|\sigma_{p_1}|$ are different from zero as they are both equal to a product of positive coefficients. Hence $k_{uv}\not=0$. The same argument can also be used in the case when $p_1$ and $p_2$  both have even length (i.e., $|V|=2k$ for some  even integer $k$) to deduce that $k_{uv}\not=0$. Hence $u$ and $v$ are adjacent in the concentration graph $G$; thus $G$ is necessarily complete.
\end{rk}

We now give an example illustrating the main result in this paper (Theorem \ref{faithG0}).

\begin{example}
Consider a Gaussian random vector $\bX=(X_1,\ldots,X_8)'$ with covariance matrix $\Sigma$ and its associated covariance graph as given in Figure \ref{treeFaithG0}.
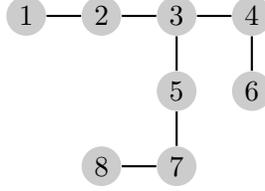
\begin{figure}[h]
\begin{center}
\begin{tikzpicture}
\tikzstyle{vertex}=[circle,fill=black!20,minimum size=15pt,inner
sep=0pt] \tikzstyle{edge} = [draw,thick,-] \foreach \pos/\name in
{{(-2,1)/1}, {(-1,1)/2},
{(0,1)/3},{(1,1)/4},{(0,0)/5},{(1,0)/6},{(0,-1)/7},{(-1,-1)/8}}
            \node[vertex] (\name) at \pos {$\name$};
    \foreach \source/ \dest  in {1/2,  2/3,3/5,3/4,4/6,5/7,7/8}
       \path[edge] (\source) -- (\dest);
\end{tikzpicture}
\caption{An $8-$vertex covariance tree $G_0$.}\label{treeFaithG0}
\end{center}
\end{figure}
Consider the sets $A=\{1,2\}$, $B=\{5\}$ and $S=\{4,6\}$. Note that $S$
does not separate $A$ and $B$ in $G_0$ as any path from $A$ and $B$ does not intersect $S$. In this case we cannot
 use the covariance global Markov property to claim that $\bX_A$ is not  independent of $\bX_B$ given $\bX_{V\setminus(A\cup B\cup S)}$.  This is because the covariance global Markov property allows us to read conditional independences present in a distribution if a separation is present in the graph. It is not an ``if and only if" property in the sense that the lack of a separation in the graph does not necessarily imply the lack of the corresponding conditional independence. We shall show however that in this example that $\bX_A$ is indeed not independent of $\bX_B$ given $\bX_{V\setminus(A\cup B\cup S)}$. In other words we shall show that the graph has the ability to capture this conditional dependence present in the probability distribution $P$.

Let us now examine the relationship between $X_2$ and $X_5$ given $\bX_{\{3,7,8\}}$. Note that in this example $V\setminus(A\cup B\cup S)=\{3,8,7\}$, $2\in A$ and $5\in B$. Note that the covariance graph associated with the probability distribution of the random vector $(X_2,X_5,\bX_{\{3,8,7\}})'$ is the subgraph represented in Figure \ref{treeFaithG0B} and can be obtained directly as a subgraph of $G_0$ induced by the subset $\{2,5,3,7,8\}$.
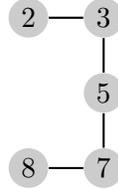
\begin{figure}[h]
\begin{center}
\begin{tikzpicture}
\tikzstyle{vertex}=[circle,fill=black!20,minimum size=15pt,inner
sep=0pt] \tikzstyle{edge} = [draw,thick,-] \foreach \pos/\name in {
{(-1,1)/2}, {(0,1)/3},{(0,0)/5},{(0,-1)/7},{(-1,-1)/8}}
            \node[vertex] (\name) at \pos {$\name$};
    \foreach \source/ \dest  in {  2/3,3/5,5/7,7/8}
       \path[edge] (\source) -- (\dest);
\end{tikzpicture}
\caption{the covariance graph
$(G_0)_{\{2,5,3,8,7\}}$}\label{treeFaithG0B}
\end{center}
\end{figure}

Since $2$ and $5$ are connected  by exactly one path in $(G_0)_{\{2,5,3,7,8\}}$, that is $p=(2,3,5)$, then the coefficient $k_{25\mid 387}$, i.e., the coefficient between $2$ and $5$ in inverse of the covariance matrix of $(X_2,X_5,\bX_{\{3,8,7\}})'$, can be computed using Theorem \ref{kgraphs:theo1} as follows
\begin{equation}\label{examptree}
k_{25\mid 387}=(-1)^{2+1} \sigma_{23}\, \sigma_{35}
\displaystyle\frac{|\Sigma(\{8,7\})|}{|\Sigma(\{2,5,3,8,7\})|}
\end{equation}
where $\Sigma(\{7,8\})$ and $\Sigma(\{2,5,3,8,7\})$ are respectively
the covariance matrices of the Gaussian random vectors $(X_7,X_8)'$
and $(X_2,X_5,\bX_{\{3,8,7\}})'$. Hence $k_{25\mid
387}\not=0$ since the right hand side of the equation in
(\ref{examptree}) is different from zero. Hence $X_2\ncind X_5\mid
\bX_{\{3,8,7\}}$.

Now recall that for any Gaussian random vector vector  $\bX_V=(X_u,\, u\in V)'$ ,
\begin{equation}\label{condGausIndep}
 \bX_A\cind \bX_B\mid \bX_{C}\mbox{ if and only if } \;
\forall\, (u,v)\in A\times B,\;\; X_u\cind X_v\mid \bX_{C}
\end{equation}
where $A$, $B$ and $C$ are pairwise disjoint subsets of $V$. The contrapositive  of (\ref{condGausIndep}) yields
$$X_2\ncind X_5\mid \bX_{\{3,7,8\}}\;\Rightarrow\;\bX_{\{1,2\}}\ncind X_5\mid \bX_{\{3,7,8\}}.$$

Hence we conclude that since $\{3,7,8\}$ does not separate $\{1,2\}$ and  $\{5\}$ therefore $\bX_{\{1,2\}}$ is not independent of $X_5$ given $\bX_{\{3,7,8\}}$, i.e.,

$$ \{1,2\}\not\perp_{G_0}\{5\}\mid \{3,7,8\} \Rightarrow \bX_{\{1,2\}}\ncind X_5\mid \bX_{\{3,7,8\}}$$.
\end{example}
\vspace{-0.4cm}
We now proceed to the proof of Theorem \ref{faithG0}.
\begin{pf} of Theorem \ref{faithG0}.
Without loss of generality we assume that $G_0$ is a connected tree. Let us assume to the contrary that $P$ is
not covariance faithful to $G_0$, then there exists a
triplet $(A,B,S)$ of pairwise disjoint subsets of $V$, such that $\bX_A\cind \bX_B\mid \bX_{V\setminus(A\cup B\cup
S)}$, but $S$ does not separate $A$ and $B$ in $G_0$, i.e.,
$$\bX_A\cind
\bX_B\mid \bX_{V\setminus(A\cup B\cup S)}\;\mbox{ and }\;
A\not\perp_{G_0} B\mid S$$

As $S$ does not separate $A$ and $B$ and since $G_0$ is a connected tree, then there exists a pair of vertices $(u,v)\in
A\times B$ such that the single path $p$ connecting $u$ and $v$ in $G_0$
 does not intersect $S$, i.e., $S\cap p=\emptyset$. Hence
$p\subseteq V\setminus S=(A\cup B)\cup (V\setminus(A\cup B\cup S))$.
Thus two cases are possible with regards to where the path $p$ can lie : either $p\subseteq A\cup B$ or
$p\cap(V\setminus(A\cup B\cup S))\not=\emptyset$.  Let us examine both cases separately.

\begin{itemize}

\item  \textbf{Case 1 :} $p\subseteq A\cup B$

In this case the entire path between $u$ and $v$ lies in $A\cup B$ and hence we can find a pair of vertices\footnote{As an illustration of this point consider the graph presented in Figure \ref{treeFaithG0}. Let  $A=\{1,2\}$, $B=\{3,5\}$ and $S=\{4,6\}$. We note that the path $p=(1,2,3,5)$ lies entirely in $A\cup B$ and hence we can find two vertices, namely, $2\in A$ and $3\in B$, belonging to path $p$ that are adjacent in $G_0$.} $(u',v')$ belonging to $p$ and $(u',v')\in A\times B$ such that $u'\sim_{G_0} v'$.\\

Recall that since $G_0$ is a tree, any induced graph of $G_0$ by a subset
of $V$ is a union of tree connected components (see Lemma \ref{TreeCon}). Hence the subgraph
$(G_0)_{W}$ of $G_0$ induced by  $W=\{u',v'\}\cup V\setminus(A\cup
B\cup S)$ is a union of tree connected components. As $u'$ and $v'$
are adjacent in $G_0$, they are also adjacent in  $(G_0)_{W}$ and
 belong to the same connected component\footnote{In our example in Figure \ref{treeFaithG0} with $W=\{2,3,8,7\}$, $(G_0)_W$ consists a union of two connected components with its respective vertices being  $\{2,3\}$ and $\{8,7\}$.} of $(G_0)_{W}$. Hence the only path between $u'$ and $v'$ is precisely the edge $(u',v')$. Using theorem  \ref{kgraphs:theo1} to compute the coefficient $k_{u'v'\mid V\setminus(A\cup B\cup S)}$, i.e.,
$(u',v')th$ coefficient in the inverse of the covariance matrix of the
random vector $\bX_W=(X_w,\,w\in W)'=(X_{u'},X_{v'},\bX_{V\setminus(A\cup B\cup S)})'$, we
obtain,
\begin{equation}\label{proofaithG0eq2a}
k_{u'v'\mid V\setminus(A\cup B\cup S)}=(-1)^{1+1} \sigma_{u'v'}\,
\displaystyle\frac{\left|\Sigma(W\setminus\{u',v'\})\right|}{|\Sigma(W)|},
\end{equation}
where $\Sigma(W)$ denotes the covariance matrix of $\bX_W$, and
$\Sigma(W\setminus\{u',v'\})$ denotes the matrix $\Sigma(W)$ with the rows and
the columns corresponding to variables $X_{u'}$ and $X_{v'}$
omitted. We can therefore deduce from (\ref{proofaithG0eq2a}) that
$k_{u'v'\mid V\setminus(A\cup B\cup S)}\not=0$. Recall that at the start of the proof we assumed to the contrary that $\bX_A\cind \bX_B\mid \bX_{V\setminus(A\cup B\cup S)}$. Now since $P$ is Gaussian, for pairwise disjoint subsets $A, B, V\setminus(A\cup B\cup C)$ then
\begin{equation}\label{proofaithG0eq1}
\bX_A\cind \bX_B\mid \bX_{V\setminus(A\cup B\cup C)} \Leftrightarrow \forall \, (u,v)\in A\times B,\; X_u\cind X_v\mid \bX_{V\setminus(A\cup B\cup C)}
\end{equation}
Note however that we have established that $X_{u'}\ncind X_{v'}\mid \bX_{V\setminus(A\cup B\cup S)}$ since $k_{u'v'\mid V\setminus(A\cup B\cup S)} \ne 0$. Hence we obtain a contradiction to (\ref{proofaithG0eq1}) since $u'\in
A$ and $v'\in B$.\\

\item  \textbf{Case 2 :} $p\cap(V\setminus(A\cup B\cup
S))\not=\emptyset$ \& $V\setminus(A\cup B\cup S)$ is not empty.

Now if $V\setminus(A\cup B\cup S)$ is empty then
$p$ has to lie entirely in $A\cup B$. This is because by assumption $p$ does not intersect $S$. The case when $p$ lies in $A\cup B$ is covered in Case 1 and hence it is assumed that $V\setminus(A\cup B\cup S)\not=\emptyset.$ \footnote{As an illustration of this point consider once more the graph presented in Figure \ref{treeFaithG0}. Consider $A=\{1,2\}$, $B=\{7,8\}$ and $S=\{4,6\}$. Here $V\setminus(A\cup B\cup S)=\{3,5\}$ and the path $p=(1,2,3,5,7,8)$ connecting $A$ and $B$ intersects $V\setminus(A\cup B\cup S)$.}

In this case there exists  a pair of vertices $(u',v')\in A\times B$ with $u',v'\in p$, such that  the vertices $u'$ and $v'$ are connected by exactly one path $p'\subseteq p$ in the induced graph $(G_0)_W$ of $G_0$ by $W=\{u',v'\}\cup V\setminus(A\cup B\cup S)$ (see Lemma \ref{TreeCon}) \footnote{In our example in figure 1 with $A=\{1,2\}$, $B=\{7,8\}$ and $S=\{4,6\}$ , the vertices $u'$ and $v'$ will correspond  to vertices $2$ and $7$ respectively, and $p'=(2,3,5,7)$, which is a path entirely contained in $V\setminus(A\cup B\cup S) \cup \{u', v'\}$.}.

Let us now use Theorem  \ref{kgraphs:theo1} to compute the coefficient
$k_{u'v'\mid V\setminus(A\cup B\cup S)}$, i.e., the
$(u',v')-$coefficient in the inverse of the covariance matrix of the
random vector $\bX_W=(X_w,\,w\in
W)'=(X_{u'},X_{v'},\bX_{V\setminus(A\cup B\cup S)})'$. We obtain that

\begin{equation}\label{proofaithG0eq2}
k_{u'v'\mid V\setminus(A\cup B\cup S)}=(-1)^{|p'|+1} |\sigma_{p'}|\,
\displaystyle\frac{\left|\Sigma(W\setminus p')\right|}{|\Sigma(W)|},
\end{equation}

where $\Sigma(W)$ denotes the covariance matrix of $\bX_W$ and
$\Sigma(W\setminus p')$ denotes $\Sigma(W)$ with the rows and the
columns corresponding to variables in path $p'$ omitted. One can therefore
 easily deduce from (\ref{proofaithG0eq2}) that $k_{u'v'\mid
V\setminus(A\cup B\cup S)}\not=0$. Thus $X_{u'}$ is not independent
of $X_{v'}$ given  $\bX_{V\setminus(A\cup B\cup S)}$. Hence once more
 we obtain a contradiction to (\ref{proofaithG0eq1}) since $u'\in
A$ and $v'\in B$.
\end{itemize}
\end{pf}

\vspace{-0.4cm}
\begin{rk} The dual result of the theorem above for the case of concentration trees was proved by \cite{Becker2005}. We note however that the argument used in the proof of Theorem \ref{faithG0} cannot also be used to prove faithfulness of Gaussian distributions that have trees as concentration graphs. The reason for this is as follows. In our proof we employed the fact that the sub-graph $(G_0)_{\{u,v\}\cup S}$ of $G_0$ induced by a subset ${\{u,v\}\cup S}\subseteq V$ is also the covariance graph associated with the Gaussian sub-random vector of $\bX_V$ as denoted by $\bX_{\{u,v\}\cup S}=(X_w,\, w\in \{u,v\}\cup S)'$. Hence it was possible to compute the coefficient $k_{uv\mid S}$ which quantifies the conditional (in)dependence between $u$ and $v$ given $S$, in terms of the paths in $(G_0)_{\{u,v\}\cup S}$ and the coefficients of the covariance  matrix of $\bX_{\{u,v\}\cup S}=(X_w,\, u\in \{u,v\}\cup S)'$. On the contrary, in the case of concentration graphs the sub-graph $G_{\{u,v\}\cup S}$ of the concentration graph $G$ induced by $\{u,v\}\cup S$ is not in general the concentration graph of the random vector $\bX_{\{u,v\}\cup S}=(X_w,\, u\in \{u,v\}\cup S)'$. Hence our approach is not directly applicable in the concentration graph setting.
\end{rk}

\section{Conclusion}

Faithfulness of a probability distribution to a graph is a crucial
assumption that is often made in the probabilistic treatment of
graphical models. This assumption describes the ability of a
graph to reflect or encode the multivariate dependencies that are present in a
joint probability distribution. Much of the methodology in this area
often do not undertake a detailed analysis of the faithfulness
assumption, as such an endeavor requires a more careful and rigorous probabilistic
study of the joint distribution at hand. In this note we looked at the
class of multivariate Gaussian distributions that are Markov with
respect to covariance graphs and prove that Gaussian distributions
which have trees as their covariance graphs are necessarily faithful.
The method of proof  that is employed in this paper is novel in the
sense that it is self contained and yields a completely new approach to demonstrating
faithfulness - as compared to the methods that are traditionally
used in the literature. Moreover, it is also vastly different in nature
from the proof of the analogous result for concentration graph models.
Hence the approach used in this paper promises to have further implications
and give other insights. Future research in this area will explore if the techniques used in this paper can be modified to prove or disprove faithfulness for other classes of graphs.

\section*{Acknowledgments}
The authors gratefully acknowledge the faculty at Stanford University for their feedback and tremendous enthusiasm for this work.

\end{document}